\documentclass{article}
\usepackage{authblk}
\usepackage{geometry}
\usepackage[titletoc]{appendix}
\usepackage{enumitem}
\usepackage{booktabs,diagbox}
\usepackage{hyperref}
\usepackage{xcolor}

\usepackage{tikz} 
\usetikzlibrary{cd}
\usetikzlibrary{arrows}

\usepackage{amssymb,amsmath,amsthm,mathtools,slashed,bm} 
\def\-{\raisebox{.75pt}{-}}

\numberwithin{table}{section}
\numberwithin{equation}{section}

\theoremstyle{definition}
\newtheorem{defn}{Definition}[section]

\newtheorem{rmk}{Remark}[section]
\theoremstyle{plain}
\newtheorem{lem}{Lemma}[section]
\newtheorem{prop}{Proposition}[section]
\newtheorem{thm}{Theorem}[section]
\newtheorem{cor}{Corollary}[section]
\newenvironment{claim}[1]{\par\noindent\underline{Claim:}\space#1}{}
\newenvironment{claimproof}[1]{\par\noindent\underline{Proof:}\space#1}{\hfill $\blacksquare$}

\usepackage{amsrefs}

\title{Monoidally Graded Manifolds}
\author[1]{Shuhan Jiang}
\date{}
\affil[1]{Max Planck Institute for Mathematics in the Sciences, 04103 Leipzig}

\begin{document}

\maketitle

\begin{abstract}
	We give a generalization of the theory of $\mathbb{Z}_2$-graded manifolds to a theory of $\mathcal{I}$-graded manifolds, where $\mathcal{I}$ is a commutative semi-ring with some additional properties. We prove Batchelor's theorem in this generalized setting. To our knowledge, such a proof is still missing except for some special cases.
\end{abstract}

\section{Introduction}

The notion of a supermanifold appeared in the late 70s when mathematicians tried to understand the concept of supersymmetry proposed by physicists \cite{Kostant77}. It extends the notion of a manifold $M$ naturally by attaching Grassmann algebras locally to $M$. The mysterious anticommutativity property of a fermionic field over $M$ can be then interpreted in terms of the anticommutativity of Grassmann algebras. When multiplying two fermionic fields, one gets a bosonic field. This process can be tracked by assigning $0 \in \mathbb{Z}_2$ to bosonic fields and $1 \in \mathbb{Z}_2$ to fermionic fields. Hence, supermanifolds are also called $\mathbb{Z}_2$-graded manifolds. Though the grading in this case is merely used to distinguish commutative and anticommutative objects.

In 1982, Witten published a seminal paper relating Morse theory to supersymmetric quantum mechanics \cite{Witten82}. It was realized since then that there exists a very deep connection between supersymmetric theories in physics and cohomology theories in mathematics. To establish such a connection, one needs to update the language of $\mathbb{Z}_2$-graded manifolds to the language of $\mathbb{Z}$-graded (or graded) manifolds \cite{Cattaneo11,Fairon17}. One major achievement in that direction is the AKSZ formalism of topological quantum field theories \cite{Alexandrov1997}, where the topological sigma models \cite{Witten88b} are reinterpreted in the language of so-called $Q$-manifolds.\footnote{A $Q$-manifold is a graded manifold equipped with a vector field $Q$ of degree $1$ satisfying $Q^2=0$.}

In cohomological field theories (or topological quantum field theories of Witten type), one can obtain useful invariants of smooth manifolds by studying observables $\mathcal{O}^{(p)}$ satisfying the following descent equations \cite{Witten88b,Witten88a}
\begin{align}\label{desceq}
	Q(\mathcal{O}^{(p)})=d(\mathcal{O}^{(p-1)}),
\end{align}
for $p \geq 1$ with $Q\mathcal{O}^{(0)}=0$, where $d$ is the de Rham differential. (\ref{desceq}) is equivalent to saying that $\mathcal{O}:=\sum_p \mathcal{O}^{(p)}$ is closed in the total complex of some bicomplex with horizontal differential $d$ and vertical differential $Q$. Such a bicomplex can be obtained by applying a ``change of coordinates'' to the variational bicomplex of a fiber bundle \cite{Jiang22}. It is then also interesting to study $\mathbb{Z}\times\mathbb{Z}$-graded (or bigraded) manifolds. 

In this paper, we follow the algebraic-geometric approaches in \cite{Kostant77,Leites80,Manin97,Carmeli11,Bartocci12,Kessler19} to give a definition of $\mathcal{I}$-graded manifolds, where $\mathcal{I}$ is an arbitrary commutative semi-ring with some additional properties. We also apply the techniques in \cite{Manin97} to give a proof of Batchelor's theorem, i.e., that every $\mathcal{I}$-graded manifold can be obtained from an $\mathcal{I}$-graded vector bundle. To our knowledge, such a proof is still missing except for some special cases \cite{Batchelor79,Covolo16,Kotov21}.
 
\section{Commutative Monoids and Parity Functions}

Let $(\mathcal{I},0,+)$ be a commutative monoid. Let $\mathbb{Z}_q$ denote the cyclic group of order $q$.
\begin{defn}
	A parity function is a (non-trivial) monoid homomorphism $p: \mathcal{I} \rightarrow \mathbb{Z}_2$.
\end{defn}
Not every $\mathcal{I}$ has a non-trivial parity function. For example, there is no non-trivial homomorphism from $\mathbb{Z}_q$ to $\mathbb{Z}_2$ when $q$ is a odd. Let $\mathcal{I}_a$ denote $p^{-1}(a)$ for $a \in \mathbb{Z}_2$.\footnote{We say that $\mathcal{I}_0$ is the even part of $\mathcal{I}$, and that $\mathcal{I}_1$ is the odd part of $\mathcal{I}$. We also say that an element of $\mathcal{I}_a$ has parity $a$ for $a=0,1$.} We have $\mathcal{I}_a + \mathcal{I}_b \subseteq \mathcal{I}_{a+b}$. Recall that an element $x$ in $\mathcal{I}$ is called cancellative if $x+y = x+z$ implies $y=z$ for all $y$ and $z$ in $\mathcal{I}$. Suppose that there is a cancellative element in $\mathcal{I}_1$. It is easy to see that such an element induces an injective map from $\mathcal{I}_{a}$ to $\mathcal{I}_{a+1}$. It follows from the Cantor-Bernstein theorem that there exists a bijection between $\mathcal{I}_0$ and $\mathcal{I}_1$. A monoid is called cancellative if every element in it is cancellative. We have shown that
\begin{prop}\label{ccm}
	Let $\mathcal{I}$ be an commutative cancellative monoid. If $\mathcal{I}$ has a non-trivial parity function $p$, then the submonoid $\mathcal{I}_0$ and its complement $\mathcal{I}_1$ have the same cardinality.
\end{prop}
\begin{rmk}
	In the finite case, proposition \ref{ccm} is no longer true if we drop the cancellative condition. For example, we can consider the commutative monoid defined by the following table.
	\begin{table}[!htbp]
		\centering
		\begin{tabular}{c|c|c|c}
			\hline
			& $0$ & $a$ & $b$ \\
			\hline
			$0$ & $0$ & $a$ & $b$ \\
			\hline
			$a$ & $a$ & $b$ & $a$ \\
			\hline
			$b$ & $b$ & $a$ & $b$ \\
			\hline
		\end{tabular}    
		\caption{A commutative non-cancellative monoid of order $3$.}\label{cncm5}
	\end{table}
	A non-trivial $p$ is defined by setting $p(0)=p(b)=0$ and $p(a)=1$.
\end{rmk}
The question now is, given an appropriate commutative cancellative monoid $\mathcal{I}$, how can one construct a parity function for it? If $\mathcal{I}$ is a finite, it is not hard to show that $\mathcal{I}$ is actually an abelian group. The fundamental theorem of finite abelian groups then tells us that $\mathcal{I}$ is isomorphic to a direct product of cyclic groups of prime-power order. By Proposition \ref{ccm}, one of these cyclic groups must be $\mathbb{Z}_{2^k}$, $k \geq 1$. We can write
$
	\mathcal{I} = \mathbb{Z}_{2^k} \times \cdots
$
and define $p$ by sending $(x, \cdots) \in \mathcal{I}$ to $a - 1 \pmod 2$, where $a$ is the order of $x \in \mathbb{Z}_{2^k}$. If $\mathcal{I}$ is infinite, the construction of $p$ is hard, perhaps not possible in general. However, one can easily work out the case when $\mathcal{I}$ is free. ($\mathcal{I}$ is then cancellative, but not a group.) Let $\mathcal{I}_0$ be the submonoid of elements generated by even number of generators. Let $\mathcal{I}_1$ be the subset of elements generated by odd number of generators. Note that $\mathcal{I}_a + \mathcal{I}_b \subseteq \mathcal{I}_{a+b}$. We obtain a parity function which sends elements in $\mathcal{I}_a$ to $a$. As an example, let $\mathcal{I}$ be $\mathbb{N}$, the monoid of natural numbers under addition. $p$ is then defined by sending even numbers to $0$ and odd numbers to $1$. 

Let $K(\mathcal{I})$ denote the Grothendieck group of $\mathcal{I}$. Recall that it can be constructed as follows. Let $\sim$ be the equivalence relation on $\mathcal{I} \times \mathcal{I}$ defined by $(a_1,a_2) \sim (b_1,b_2)$ if there exists a $c \in \mathcal{I}$ such that $a_1 + b_2 + c = a_2 + b_1 + c$. The quotient $K(\mathcal{I}) = \mathcal{I} \times \mathcal{I} / \sim$ has a group structure by $[(a_1,a_2)] + [(b_1,b_2)]=[(a_1+b_1,a_2+b_2)]$.
\begin{prop}\label{gropar}
	Let $p$ be a parity function for $\mathcal{I}$. The map
	\begin{align*}
		p': K(\mathcal{I}) &\rightarrow \mathbb{Z}_2 \\
		[(a_1,a_2)] &\mapsto p(a_1) + p(a_2)
	\end{align*}
	is well-defined and gives a parity function for $K(\mathcal{I})$.
\end{prop}
\begin{rmk}
	When $\mathcal{I}$ is cancellative, it can be seen as a submonoid of $K(\mathcal{I})$ by the embedding
	\begin{align*}
		\iota: \mathcal{I} &\rightarrow K(\mathcal{I}) \\
		a &\mapsto [(a,0)].
	\end{align*}
	For this reason, we sometimes simply write $a-b$ to denote $[(a,b)] \in K(\mathcal{I})$. The cancellative property is not necessary for the proof of Proposition \ref{gropar}. But it guarantees the non-triviality of $p'$, since $p'$ restricted to $\mathcal{I}$ must coincide with $p$. 
\end{rmk}
\begin{proof}
	Let $(a_1,a_2)$ and $(b_1,b_2)$ represent the same element of $K(\mathcal{I})$, i.e., there exist some $c$ such that $a_1 + b_2 + c = a_2 + b_1 + c$. One then concludes that $a_1+b_2$ and $a_2+b_1$ must have the same parity. Note that, for $a,b \in \mathbb{Z}_2$, $a=b$ if and only if $a+b=0$. We have
	\begin{align*}
		p'([(a_1,a_2)]) + p'(([b_1,b_2)]) = p(a_1+b_2) + p(a_2+b_1) = 0.
	\end{align*}
	Hence $p'([a_1,a_2])=p'([b_1,b_2])$.
\end{proof}
As an example, consider $K(\mathbb{N}) = \mathbb{Z}$, the monoid of integers under addition. The parity function $p'$ induced from the parity function $p$ for $\mathbb{N}$ again sends even numbers to $0$ and odd numbers to $1$. 

\section{Monoidally Graded Ringed Spaces}

Let $R$ be a commutative ring. Let $\mathcal{I}$ be a countable commutative cancellative monoid equipped with a parity function $p$.
\begin{defn}
	An $\mathcal{I}$-graded $R$-module is an $R$-module $V$ with a family of sub-modules $\{V_i\}_{i \in \mathcal{I}}$ indexed by $\mathcal{I}$ such that $V=\bigoplus_{i \in \mathcal{I}} V_i$. $v \in V$ is said to be homogeneous if $v \in V_i$ for some $i \in \mathcal{I}$. We use $d(v)$ to denote the degree of $v$, $d(v)=i$. 
\end{defn} 
Given two $\mathcal{I}$-graded $R$-modules $V$ and $W$, we make the direct sum $V \oplus W$ and the tensor product $V \otimes W$ into $\mathcal{I}$-graded $R$-modules by setting
\begin{align*}
	V \oplus W = \bigoplus_{i \in \mathcal{I}}(V_i \oplus W_i), \quad V \otimes W = \bigoplus_{k \in \mathcal{I}}\left(\bigoplus_{i+j=k} V_i \otimes W_j\right).
\end{align*}
We can also make the space $\mathrm{Hom}(V,W)$ of $R$-linear maps from $V$ to $W$ into a $K(\mathcal{I})$-graded $R$-module by setting
\begin{align*}
	\mathrm{Hom}(V,W) = \bigoplus_{\alpha \in K(\mathcal{I})} \mathrm{Hom}(V,W)_{\alpha}, \quad \mathrm{Hom}(V,W)_{\alpha}=\{f \in \mathrm{Hom}(V,W)|f(V_{i}) \subset W_{j}, [(j,i)]=\alpha\}.
\end{align*}
A morphism from $V$ to $W$ is just an element of $\mathrm{Hom}(V,W)_{0}$, i.e., an $R$-linear map of degree $0$. 
\begin{rmk}
	$\mathrm{Hom}(V,W)$ is in general not $\mathcal{I}$-graded. This is because that we should assign degree ``$j-i$'' to a map $f$ which maps elements in $V_i$ to elements in $w \in W_j$. But the minus operation does make sense for a general monoid $\mathcal{I}$. So we have to work with $K(\mathcal{I})$, the group completion of $\mathcal{I}$. Note that $V^*=\mathrm{Hom}(V,R)$, the dual of $V$, is in particular $K(\mathcal{I})$-graded. (The degree of elements in $V_i^*$ is $-i$.) Hence $V^* \otimes W$, which is isomorphic to $\mathrm{Hom}(V,W)$, is $K(\mathcal{I})$-graded by assigning degree $j-i$ to elements in $V_i^* \otimes W_j$. Everything is consistent.
\end{rmk}
Now, suppose that $\mathcal{I}$ also has a commutative multiplicative structure which is compatible with the additive structure. That is, it is a commutative cancellative semi-ring. We write $ab$ as the multiplication of $a$ and $b$ in $\mathcal{I}$.
\begin{defn}
	An $\mathcal{I}$-graded $R$-module $A$ is called an $\mathcal{I}$-graded $R$-algebra if $A$ is a unital associative $R$-algebra and if the multiplication $\mu: A \otimes A \rightarrow A$ is a morphism of $\mathcal{I}$-graded $R$-modules. We write $xy=\mu(x\otimes y)$ as the shorthand notation for multiplications of $A$. $A$ is said to be commutative if
	\begin{align}\label{grasign}
		xy - (-1)^{p(x)p(y)} yx = 0
	\end{align}
	for all homogeneous $x,y \in A$, where we use $p(x)p(y)$ to denote $p(d(x)d(y)) \in \mathbb{Z}_2$.
\end{defn}
\begin{rmk}
	Here we have to be careful about the sign factor appearing in the right hand side of (\ref{grasign}). Although both of $\mathcal{I}$ and $\mathbb{Z}_2$ are semi-rings\footnote{The multiplicative structure on $\mathbb{Z}_2$ is inherited from the one on $\mathbb{Z}$.}, $p$ is not necessarily a semi-ring homomorphism and we do not have $p(d(x)d(y)) = p(d(x))p(d(y))$ in general. To choose which as the sign factor is just a matter of convention.
\end{rmk}
Morphisms of $\mathcal{I}$-graded algebras are simply linear maps of degree $0$ which preserves the algebraic structures. We use $\mathrm{Comm}$-$\mathrm{Alg}_{\mathcal{I}}$ to denote the category of commutative $\mathcal{I}$-graded algebras.
\begin{defn}
	The tensor algebra $\mathrm{T}(V)$ is the $\mathcal{I}$-graded $R$-module $\mathrm{T}(V)=\bigoplus_{n \in \mathbb{N}} V^{\otimes^n}$, together with the tensor product $\otimes$ as the canonical multiplication. The symmetric algebra $\mathrm{S}(V)$ is the quotient algebra of $\mathrm{T}(V)$ by the $\mathcal{I}$-graded two-sided ideal generated by
	\begin{align*}
		v \otimes w - (-1)^{p(v)p(w)}w \otimes v,
	\end{align*}
	where $v, w \in V \subset \mathrm{T}(V)$ are homogeneous.
\end{defn}
\begin{rmk}
	$\mathrm{S}(V)$ has a canonical $\mathbb{N}$-grading inherited from $\mathrm{T}(V)$ which should not be confused with its $\mathcal{I}$-grading. We write $\mathrm{S}(V) = \bigoplus_{n \in \mathbb{N}} \mathrm{S}^n(V)$ to indicate that fact. Note that $\mathrm{S}^0(V) = R$, but $\mathrm{S}(V)_0$, the sub-space of homogeneous elements of degree $0$, is in general larger than $R$.
\end{rmk}
$\mathrm{S}(V)$ is universal in the sense that, given a commutative $\mathcal{I}$-graded $R$-algebra $A$ and a morphism $f: V \rightarrow A$. There exists a unique algebraic homomorphism $\tilde{f}: \mathrm{S}(V) \rightarrow A$ such that the following diagram commutes
\[
\begin{tikzcd}
	V \arrow{dr}{f} \arrow[hook]{r}{\iota} & \mathrm{S}(V) \arrow{d}{\tilde{f}} \\
	& A
\end{tikzcd}
\]
where $\iota: V \rightarrow \mathrm{S}(V)$ is the canonical embedding. Note that $\tilde{f}$ preserves the $\mathcal{I}$-grading, i.e., it is a morphism in $\mathrm{Comm}$-$\mathrm{Alg}_{\mathcal{I}}$. Choosing $A$ to be $R$ (viewed as an $\mathcal{I}$-graded $R$-algebra whose components of non-zero degree are $0$.) and $f$ to be the zero map, we obtain an $R$-algebra homomorphism from $\mathrm{S}(V)$ to $R$. We denote this map by $\epsilon$. Note that $\ker \epsilon = \bigoplus_{n>0} \mathrm{S}^n(V)$.

Let $k$ be a field and $R$ be a commutative $k$-algebra. Let $A$ be a commutative $\mathcal{I}$-graded $k$-algebra.
\begin{defn}
	A $k$-algebra epimorphism $\epsilon: A \rightarrow R$ is called a body map of $A$ if $\ker \epsilon \supset I$, where $I$ is the ideal in $A$ generated by homogeneous elements of non-zero degree.
\end{defn}
By definition, $\epsilon$ preserves the $\mathcal{I}$-grading of $A$.
\begin{defn}
	Let $\epsilon$ be a body map of $A$. $A$ is said to be projected if the short exact sequence
	\begin{align*}
		0 \longrightarrow \ker \epsilon \longrightarrow A \xrightarrow{~\epsilon~} R \longrightarrow 0
	\end{align*}
	splits.
\end{defn}
The splitting gives $A$ an $R$-module structure depending on $\epsilon$, with respect to which $\epsilon$ becomes an $R$-algebra homomorphism. Conversely, $A$ is projected if $A$ has an $R$-module structure and $\epsilon$ preserves that structure.
\begin{lem}\label{unieps}
	Let $V$ be an $\mathcal{I}$-graded $R$-module with $V_0=0$. Let $\epsilon$ be an $R$-linear body map of $\mathrm{S}(V)$. Then $\epsilon$ is unique.
\end{lem}
\begin{proof}
	In this case, $\mathrm{S}(V) = R \oplus I$ where  $I=\bigoplus_{n>0} \mathrm{S}^n(V)$. Since $I \subset \ker \epsilon$ and $\epsilon$ is $R$-linear, the only possible choice of $\epsilon$ is the canonical one.
\end{proof}	 
\begin{rmk}\label{uniepsrmk}
	Let $V$ be as in Lemma \ref{unieps}. Suppose $A \cong \mathrm{S}(V)$ as $\mathcal{I}$-graded $k$-algebras. In particular, this implies that $A$ admits a decomposition $A = A' \oplus I$ where $A' \cong R$ and $I$ is the ideal generated by homogeneous elements of non-zero degree. Let $\epsilon$ be a body map of $A$. Since $I \subset \ker \epsilon$, $\epsilon$ is determined by $\epsilon|_{A'}$. In other words, $\epsilon$ is determined by a $k$-algebra endomorphism of $R$.
\end{rmk}
More can be said if $V$ is free.
\begin{lem}\label{isosym}
	Let $V$ be a free $\mathcal{I}$-graded $R$-module with $V_0=0$. Let $\epsilon$ be an $R$-linear body map of $\mathrm{S}(V)$. (By Lemma \ref{unieps}, $\epsilon$ is the canonical one.) Let $I$ denote the kernel of $\epsilon$. Then there exists an $R$-algebra isomorphism
	\begin{align*}
		\mathrm{S}(V) \cong \mathrm{S}(I/I^2),
	\end{align*}
	where $I^2$ is the square of the ideal $I$.
\end{lem}
\begin{proof}
	Let $\iota: V \hookrightarrow \mathrm{S}(V)$ be the canonical embedding. Since $I=\bigoplus_{n>0} \mathrm{S}^n(V)$, we have $\iota(V) \subset I$, which yields another embedding $V \hookrightarrow I/I^2 \hookrightarrow  \mathrm{S}(I/I^2)$, which induces the desired isomorphic map between $\mathrm{S}(V)$ and $\mathrm{S}(I/I^2)$.    	
\end{proof}
\begin{defn}
The $\mathcal{I}$-graded algebra of formal power series on $V$ is the $R$-module
\begin{align*}
	\overline{\mathrm{S}(V)} = \prod_{n \in \mathbb{N}} \mathrm{S}^n(V)
\end{align*}
equipped with the canonical algebraic multiplication.
\end{defn}
\begin{rmk}
	As is in the case of $\mathcal{I}=\mathbb{Z}$ \cite{Fairon17}, it is actually crucial to work with $\overline{\mathrm{S}(V)}$ instead of $\mathrm{S}(V)$ when the even part of $V$ is non-trivial. The former allows us to have a coordinate description of morphisms between ``$\mathcal{I}$-graded domains'', a notion of partition of unity for ``$\mathcal{I}$-graded manifolds'', and more.
\end{rmk}
Let $I$ be the kernel of the canonical body map of $\mathrm{S}(V)$. One can equip $\mathrm{S}(V)$ with the so-called $I$-adic topology.\footnote{To each point $x$ of $\mathrm{S}(V)$ one assigns a collection of subsets $\mathcal{B}(x)=\{x+I^n\}_{x \in A, n>0}$. The $I$-adic topology is then the unique topology on $\mathrm{S}(V)$ such that $\mathcal{B}(x)$ forms a neighborhood base of $x$ for all $x$.} Moreover, one can consider the $I$-adic completion of $\mathrm{S}(V)$ which is defined as the inverse limit 
\begin{align*}
	\widehat{\mathrm{S}(V)}_I := \varprojlim \mathrm{S}(V)/I^n
\end{align*}
of the inverse system $((\mathrm{S}(V)/I^n)_{n \in \mathbb{N}}, (\pi_{m,n})_{n \leq m \in \mathbb{N}})$, where $\pi_{m,n}: \mathrm{S}(V)/I^m \rightarrow \mathrm{S}(V)/I^n$ is the canonical projection. Note that there is also a canonical projection $\mathrm{S}(V) \rightarrow \mathrm{S}(V)/I^n$ for each $n \in \mathbb{N}$. By the universal property of the inverse limit, one obtains a morphism 
\begin{align*}
	\iota_I: \mathrm{S}(V) \rightarrow \widehat{\mathrm{S}(V)}_I
\end{align*}
with kernel being $\bigcap_{n\geq0}I^n=\{0\}$. On the other hand, it is easy to see that $\mathrm{S}(V)/I^n \cong \bigoplus_{i=0}^{n-1} \mathrm{S}^i(V)$ for $n \geq 1$. It follows that there is a canonical isomorphism $\widehat{\mathrm{S}(V)}_I \cong \overline{\mathrm{S}(V)}$ under which $\iota_I$ coincides with the canonical inclusion $\mathrm{S}(V) \hookrightarrow \overline{\mathrm{S}(V)}$. 

In fact, $\mathrm{S}(V)$ can be made into a metric space such that $\overline{\mathrm{S}(V)}$ is the completion of $\mathrm{S}(V)$ with respect to the metric structure \cite{Singh11}. The metric-induced topology on $\overline{\mathrm{S}(V)}$, with a slight abuse of notation, coincides with the $I$-adic topology on $\overline{\mathrm{S}(V)}$, where $I=\prod_{n>0}\mathrm{S}^n(V)$. 

\begin{lem}\label{unifor}
	 Let $A$ be a commutative $\mathcal{I}$-graded $R$-algebra. Let $J$ be an ideal of $A$ such that $A$ is $J$-adic complete. $\overline{\mathrm{S}(V)}$ is universal in the sense that, given a morphism $f: V \rightarrow A$ such that $f(V) \subset J$, there exists a unique (continuous) algebraic homomorphism $\tilde{f}: \overline{\mathrm{S}(V)} \rightarrow A$ such that the following diagram commutes
	\[
	\begin{tikzcd}
		V \arrow{dr}{f} \arrow[hook]{r}{\iota} & \overline{\mathrm{S}(V)} \arrow{d}{\tilde{f}} \\
		& A
	\end{tikzcd}
	\]
\end{lem}
\begin{proof}
	We already know that $f$ induces a unique morphism $f':\mathrm{S}(V) \rightarrow A$ such that $f' \circ \iota = f$. By assumption, $f'$ extends naturally to a morphism $\tilde{f}: \overline{\mathrm{S}(V)} \rightarrow \widehat{A}_J \cong A$. 
	
	\begin{claim}
		$\tilde{f}$ is continuous.
	\end{claim}
    \begin{claimproof}
    	It suffices to show that $\tilde{f}^{-1}(J^m)$ is a neighborhood of $0$ for any $m \in \mathbb{N}$. By assumption, $I \subset \tilde{f}^{-1}(J)$. It follows that $I^m \subset \tilde{f}^{-1}(J)^m \subset \tilde{f}^{-1}(J^m)$.  
    \end{claimproof}

    Since $\mathrm{S}(V)$ is dense in $\overline{\mathrm{S}(V)}$ and $\tilde{f}|_{\mathrm{S}(V)}=f'$, $\tilde{f}$ is also unique.
\end{proof}
\begin{rmk}\label{isofor}
	Likewise, we have a canonical body map of $\overline{\mathrm{S}(V)}$ induced from the zero map $V \rightarrow R$. Similar results like Lemma \ref{unieps} and Lemma \ref{isosym} also hold. For example, we have
	\begin{align*}
		\overline{\mathrm{S}(V)} \cong \overline{\mathrm{S}(I/I^2)},
	\end{align*}
	where $V$ and $I$ are as in Lemma \ref{isosym}.
\end{rmk}
\begin{lem}\label{invfor}
	Let $\epsilon$ be the canonical body map of $\overline{\mathrm{S}(V)}$. Then for $f \in \overline{\mathrm{S}(V)}$, $f$ is invertible if and only if $\epsilon(f)$ is invertible.
\end{lem}
\begin{proof}
	``$\Rightarrow$'': Trivial. \\    
	``$\Leftarrow$'': Suppose $\epsilon(f) = c$ where $c \in R$ is invertible. We can write $f = c+f'$ where $f' \in \prod_{n\geq 1} \mathrm{S}^n(V)$. Note that $(f')^k \in \prod_{n \geq k} \mathrm{S}^n(V)$ for all $k>0$. We can then set the inverse of $f$ to be the formal sum $f^{-1}:= c^{-1} \sum_{k \in \mathbb{N}} (-1)^k (c^{-1}f')^k$. ($f^{-1}$ is well-defined because the formal sum restricted to each $\mathrm{S}^n(V)$ is a finite sum.)
\end{proof} 
\begin{cor}\label{localring}
	$\overline{\mathrm{S}(V)}$ is local if $R$ is local.
\end{cor}
\begin{proof}
	Choose a non-unit $f \in \overline{\mathrm{S}(V)}$. Let $c=\epsilon(f)$. By Lemma \ref{invfor}, $c$ is a non-unit. Since $R$ is local, $1-c$ is invertible. $1-f$ is then a unit by Lemma \ref{invfor}.
\end{proof}

Recall that a ringed space $(X,\mathcal{O})$ is a topological space $X$ with a sheaf of rings $\mathcal{O}$ on $X$.
\begin{defn}\label{Irs}
	An $\mathcal{I}$-graded ringed space is a ringed space  $(X,\mathcal{O})$ such that
	\begin{enumerate}
		\item $\mathcal{O}(U)$ is an $\mathcal{I}$-graded algebra for any open subset $U$ of $X$;
		\item the restriction morphism
		$
			\rho_{V,U}: \mathcal{O}(U) \rightarrow \mathcal{O}(V)
		$
	    is a morphism of $\mathcal{I}$-graded algebras.
	\end{enumerate}
   A morphism between two $\mathcal{I}$-graded ringed spaces $(X_1,\mathcal{O}_1)$ and $(X_2,\mathcal{O}_2)$ is just a morphism $\varphi=(\tilde{\varphi},\varphi^*)$ between ringed spaces such that
   $\varphi^*_U: \mathcal{O}_2(U) \rightarrow \mathcal{O}_1(\tilde{\varphi}^{-1}(U))$ preserves the $\mathcal{I}$-grading for any open subset $U$ of $X_2$.
\end{defn}

Let $(X,C)$ be a ringed space where $C(U)$ are commutative rings. One can define $\mathcal{I}$-graded $C$-modules and commutative $\mathcal{I}$-graded $C$-algebras in a similar way. In particular, the structure sheaf $\mathcal{O}$ of an $\mathcal{I}$-graded ringed space can be viewed as an $\mathcal{I}$-graded $C$-algebra if $C$ is a sub-sheaf of $\mathcal{O}$ such that $C(U)$ are homogeneous sub-algebras of degree $0$ of $\mathcal{O}(U)$.  
\begin{defn}
	Let $\mathcal{F}$ be an $\mathcal{I}$-graded $C$-module. The formal symmetric power $\overline{\mathrm{S}(\mathcal{F})}$ of $\mathcal{F}$ is the sheafification of the presheaf
	\begin{align*}
		U \rightarrow \overline{\mathrm{S}(\mathcal{F}(U))},
	\end{align*}
    where $\overline{\mathrm{S}(\mathcal{F}(U))}$ is the $\mathcal{I}$-graded algebra of formal power series on the $C(U)$-module $\mathcal{F}(U)$.
\end{defn}
By definition, $\overline{\mathrm{S}(\mathcal{F})}$ is a commutative $\mathcal{I}$-graded $C$-algebra.
\begin{lem}\label{uniforshf}
	Let $\mathcal{A}$ be a commutative $\mathcal{I}$-graded $C$-algebra. Let $\mathcal{B}$ be a sub-sheaf of $\mathcal{A}$ such that $\mathcal{A}(U)$ is $\mathcal{B}(U)$-adic complete for all open subsets $U$. $\overline{\mathrm{S}(\mathcal{F})}$ is universal in the sense that, given a morphism of $\mathcal{I}$-graded $C$-modules $F: \mathcal{F} \rightarrow \mathcal{A}$ such that $F(\mathcal{F}(U)) \subset \mathcal{B}(U)$ for all open subsets $U$, there exists a unique morphism of $\mathcal{I}$-graded $C$-algebras $\tilde{F}: \overline{\mathrm{S}(\mathcal{F})} \rightarrow \mathcal{A}$ such that the following diagram commutes
	\[
	\begin{tikzcd}
		\mathcal{F} \arrow{dr}{F} \arrow[hook]{r}{\iota} & \overline{\mathrm{S}(\mathcal{F})} \arrow{d}{\tilde{F}} \\
		& \mathcal{A}
	\end{tikzcd}
	\]
	where $\iota: \mathcal{F} \rightarrow \overline{\mathrm{S}(\mathcal{F})}$ is the canonical monomorphism.
\end{lem}
\begin{proof}
	This follows directly from the universal property of sheafification\footnote{That is, given a presheaf $\mathcal{F}$, a sheaf $\mathcal{G}$, and a presheaf morphism $F: \mathcal{F} \rightarrow \mathcal{G}$, there exists a unique sheaf morphism $\tilde{F}: \mathcal{F}^{\sharp} \rightarrow \mathcal{G}$ such that $\tilde{F} \circ \iota = F$, where $\mathcal{F}^{\sharp}$ is the sheafification of $\mathcal{F}$ and $\iota: \mathcal{F} \rightarrow \mathcal{F}^{\sharp}$ is the canonical morphism.} and the universal property of $\overline{\mathrm{S}(\mathcal{F}(U))}$ stated in Lemma \ref{unifor}.
\end{proof}
To end this section, we state the following lemma taken from \cite{Manin97}.
\begin{lem}\label{obspl}
	Let 
	\begin{align}\label{ghf}
		0 \longrightarrow \mathcal{G} \longrightarrow \mathcal{H} \longrightarrow \mathcal{F} \longrightarrow 0.
	\end{align}
    be a short exact sequence of $C$-modules where $\mathcal{F}$ and $\mathcal{G}$ are locally free $C$-modules. Then the obstruction of the existence of a splitting of (\ref{ghf}) can be represented as an element in the first sheaf cohomology group $H^1(X,\mathrm{Hom}(\mathcal{F},\mathcal{G}))$ of $\mathrm{Hom}(\mathcal{F},\mathcal{G})$.
\end{lem}

\section{Monoidally Graded Domains}

Throughout this section, $V$ is a real $\mathcal{I}$-graded vector space with $V_0=0$. The dimension of the homogeneous sub-space $V_i$ of $V$ is $m_i$. We also assume that only finitely many of $m_i$ are non-zero.
\begin{defn}
	Let $U$ be a domain of $\mathbb{R}^n$. An $\mathcal{I}$-graded domain $\mathcal{U}$ of dimension $n|(m_{i})_{i \in \mathcal{I}}$ is an $\mathcal{I}$-graded ringed space $(U,\mathcal{O})$, where $\mathcal{O}$ is the sheaf of $\overline{\mathrm{S}(V)}$-valued smooth functions. 
\end{defn}
\begin{rmk}
	$\mathcal{U}$ is a locally ringed space by Corollary \ref{localring}. 
\end{rmk}
For example, a domain $U$ with the sheaf $C^{\infty}$ of smooth functions on $U$ is an $\mathcal{I}$-graded domain of dimension $n|(0,\cdots)$, which is denoted again by $U$ for simplicity.
\begin{lem}\label{unicinf}
	Let $F: C^{\infty} \rightarrow C^{\infty}$ be an endomorphism of sheaves of commutative rings on $U$. Then $F$ must be the identity.
\end{lem}
\begin{proof}
	First, we show that $F$ is actually an endomorphism of sheaves of unital $\mathbb{R}$-algebras on $U$. It suffices to show that $F$ restricted to any open subset of $U$ sends a constant function to itself. We know this is true for $\mathbb{Q}$-valued constant functions. Now, if $F$ sends a constant function $f$ to a non-constant function $g$, then one can find two rational number $b_1$ and $b_2$ such that $g-b_1$ and $g-b_2$ are non-invertible. But then the pre-images $f-b_1$ and $f-b_2$ are non-invertible, which implies that $f$ is non-constant: a contradiction. To show that $g$ actually equals $f$, use the fact that the only field endomorphism of $\mathbb{R}$ is the identity.
	
	Let $p \in U$. $F$ induces a unital ring endomorphism $F_p$ on the stalk $C^{\infty}_p$. On the other hand, for any open neighborhood $U_p \subset U$ of $p$, the evaluation map 
	\begin{align*}
		\mathrm{ev}: C^{\infty}(U_p) &\rightarrow \mathbb{R} \\
		f &\mapsto f(p)
	\end{align*}
	induces a map $\mathrm{ev}_p: C^{\infty}_p \rightarrow \mathbb{R}$. For $f_p \in C^{\infty}_p$, it is easy to see that $f_p$ is invertible if and only if $\mathrm{ev}_p(f_p) \neq 0$. Let $c=\mathrm{ev}_p(F_p(f_p))$. $f_p-c$ is non-invertible. Hence $\mathrm{ev}_p(f_p) = c$. In other words, for any open subset $U'$ of $U$, we have $F_{U'}(f)(p) = f(p)$ for all $f \in C^{\infty}(U')$ and all $p \in U'$. This implies $F = \mathrm{id}$.
\end{proof}
A morphism between $\mathcal{I}$-graded domains is just a morphism of $\mathcal{I}$-graded locally ringed spaces. Recall that we have the canonical body map $\epsilon: C^{\infty}(U) \otimes \overline{\mathrm{S}(V)} \rightarrow C^{\infty}(U)$.
\begin{prop}\label{uniepshf}
	There exists a unique monomorphism $\varphi: U \rightarrow \mathcal{U}$ with $\tilde{\varphi} = \mathrm{id}$.
\end{prop}
\begin{proof}
	Existence is guaranteed by $\epsilon$. Uniqueness follows from Remark \ref{uniepsrmk} and Lemma \ref{unicinf}.
\end{proof}
We also have a canonical morphism for the other direction $\mathcal{U} \rightarrow U$ induced by the canonical embedding $\iota: C^{\infty}(U) \rightarrow C^{\infty}(U) \otimes \overline{\mathrm{S}(V)}$.\footnote{There will be no longer such a canonical morphism if we go the category of $\mathcal{I}$-graded manifolds.} Note that $\epsilon \circ \iota = \mathrm{id}$ on $C^{\infty}(U)$. 
\begin{prop}\label{epicom}
	Let $\varphi=(\tilde{\varphi},\varphi^*)$ be a morphism from $\mathcal{U}_1=(U_1,\mathcal{O}_1)$ to $\mathcal{U}_2=(U_2,\mathcal{O}_2)$. The following diagram commutes.
	\[
	\begin{tikzcd}
		\mathcal{U}_1 \arrow{r}{\varphi}  &\mathcal{U}_2  \\
		U_1 \arrow[hook]{u} \arrow{r}{\tilde{\varphi}} &U_2 \arrow[hook]{u}
	\end{tikzcd}
	\]
\end{prop}
\begin{proof}
	Let $U$ be an open subset of $U_2$. Let $f \in \mathcal{O}_2(U)$. We need to show that
	\begin{align*}
		\epsilon(\varphi^*(f)) = \epsilon(f) \circ \tilde{\varphi}.
	\end{align*}
	Suppose this does not hold. One can find a $p \in \tilde{\varphi}^{-1}(U)$ such that $\epsilon(\varphi^*(f))(p)=c\neq \epsilon(f)(\tilde{\varphi}(p))$. Then there exists an open neighborhood $U'\subset U$ of $\tilde{\varphi}(p)$ such that $\epsilon(f)-c$ is invertible. By Lemma \ref{invfor}, $f-c$ is also invertible on $U'$, which implies that $\varphi^*(f-c)$ is invertible on $\tilde{\varphi}^{-1}(U') \subset \tilde{\varphi}^{-1}(U)$, which contradicts the fact that $\epsilon(\varphi^*(f-c))$ is non-invertible on $\tilde{\varphi}^{-1}(U')$.
\end{proof}
\begin{defn}
A coordinate system of $\mathcal{U}$ is a collection of functions $(x^{\mu},\theta_{i,a})$ such that
\begin{enumerate}
	\item $x^{\mu}$ are elements of $\mathcal{O}(U)_0$ such that $\epsilon(x^{\mu})$ form a coordinate system of $U$;
	\item $\theta_{i,a}$ are homogeneous elements of $\mathcal{O}(U)$ of degree $d(\theta_{i,a})=i$, $i \neq 0$ and $a=1,\cdots,m_i$, which generate $\mathcal{O}(U)$ as a $C^{\infty}(U)$-algebra.
\end{enumerate}
\end{defn}
Suppose that $\mathcal{I}$ can be given a total order $<$. It follows that any function $f \in \mathcal{O}(U)$ can be written uniquely in the form
\begin{align}\label{decomp}
	f=\sum_{\mathcal{J}} \sum_{\beta} f_{\mathcal{J},\beta}(x^{\mu}) \prod_{j \in \mathcal{J}} \theta_j^{\beta^j},
\end{align}
where
\begin{itemize}
	\item $\mathcal{J} \in \mathrm{Pow}(\mathcal{I})$, $\beta=(\beta^j)_{j \in \mathcal{J}}$, $\beta^j=(\beta^j_1,\dots,\beta^j_{m_j})$, $\beta^j_k \in \{0,1\}$ if $p(j)=1$, $\beta^j_k \in \mathbb{N}$ if $p(j)=0$;
	\item $\theta_j^{\beta^j} = \theta_{j,1}^{\beta^j_1} \cdots \theta_{j,m_j}^{\beta^j_{m_j}}$, the product $\prod_{j \in \mathcal{J}} \theta_j^{\beta^j}$ is arranged in a proper order such that $\theta_j^{\beta^j}$ is on the left of $\theta_{j'}^{\beta^{j'}}$ whenever $j < j'$;
	\item For a smooth function $g \in C^{\infty}(U)$, the notation $g(x^{\mu})$ should be understood as
	\begin{align}\label{grasscont}
		g(x^{\mu}):=\sum_{i_1=0}^{\infty}\cdots\sum_{i_n=0}^{\infty}\frac{1}{i_1!\cdots i_n!}\partial_1^{i_1}\cdots \partial_n^{i_n}g(\epsilon(x^{\mu})) (x^1-\epsilon(x^1))^{i_1}\cdots(x^n-\epsilon(x^n))^{i_n}.
	\end{align}
	Hence, $g(x^{\mu})$ is an element in $\mathcal{O}(U)_0$ instead of $C^{\infty}(U)$.
\end{itemize}
The sum in (\ref{decomp}) is well-defined because, by assumption, only finitely many of $m_j$ are non-zero.
\begin{rmk}
	One may wonder how we obtain (\ref{decomp}). In fact, by definition, every function $f$ can be expressed in the form
	\begin{align*}
		f=\sum_{\mathcal{J}} \sum_{\beta} f_{\mathcal{J},\beta}(\epsilon(x^{\mu})) \prod_{j \in \mathcal{J}} \theta_j^{\beta^j}.
	\end{align*}
	One can then define a map from $\mathcal{O}(U)$ to itself by sending $g(\epsilon(x^{\mu}))$ to $g(x^{\mu})$. Now consider another map which sends $g(\epsilon(x^{\mu}))$ to $g^-(x^{\mu})$, where
	\begin{align*}
		g^-(x^{\mu}):=\sum_{i_1=0}^{\infty}\cdots\sum_{i_n=0}^{\infty}\frac{1}{i_1!\cdots i_n!}\partial_1^{i_1}\cdots \partial_n^{i_n}g(\epsilon(x^{\mu})) (\epsilon(x^1)-x^1)^{i_1}\cdots(\epsilon(x^n)-x^n)^{i_n}.
	\end{align*}
	Using the binomial theorem, it is easy to see that the second map is the inverse of the first. In fact, the reader may notice that the map $g(\epsilon(x^{\mu})) \mapsto g(x^{\mu})$ is exactly the ``Grassmann analytic continuation map'' defined in \cite{Rogers07}.
\end{rmk}
\begin{cor}\label{undermap}
	Let $\varphi=(\tilde{\varphi},\varphi^*)$ be as in Proposition \ref{epicom}. $\tilde{\varphi}$ is uniquely determined by $\varphi^*$.
\end{cor}
\begin{proof}
	Let $(x^{\mu},\theta_{i,a})$ be a coordinate system of $\mathcal{U}_2$. By Proposition \ref{epicom}, one has $\tilde{\varphi}^{\mu}=\epsilon(\varphi^*x^{\mu})$, where $(\tilde{\varphi}^{\mu})$ is $\tilde{\varphi}$ expressed in the coordinate system $(\epsilon(x^{\mu}))$ of $U_2$.
\end{proof}
Let $\mathrm{ev}$ be the evaluation map of $C^{\infty}(U)$ at $p \in U$. Let $s_p$ denote $\mathrm{ev} \circ \epsilon$. Let $I_p$ denote the kernel of $s_p$. We follow \cite{Leites80} to prove the following lemmas.
\begin{lem}\label{hadamard}
	For any functions $f \in \mathcal{O}(U)$ and any integer $k \geq 0$, there is a polynomial $P_k$ in the coordinates $(x^{\mu},\theta_{i,a})$ such that $f - P_k \in I^{k+1}_p$.
\end{lem}
\begin{proof}
	Use the classical Hadamard lemma and the decomposition (\ref{decomp}).
\end{proof}
\begin{lem}\label{critequal}
	Let $f$ and $g$ be functions of $\mathcal{O}(U)$, then $f=g$ if and only if $f-g \in I^k_p$ for all $k \in \mathbb{N}$ and $p \in U$. In other words, $\bigcap_{p \in U} \bigcap_{k \in \mathbb{N}} I^k_p = \{0\}$.
\end{lem}
\begin{proof}
	Let $h=f-g$. Apply the decomposition (\ref{decomp}) to $h$, then by Lemma \ref{hadamard}, $h_{\mathcal{J},\beta}=0$ for all $\mathcal{J}$ and $\beta$. Hence $h=0$.
\end{proof}
\begin{lem}\label{recovtop}
	Any morphism of $\mathcal{I}$-graded $\mathbb{R}$-algebras $s: \mathcal{O}(U) \rightarrow \mathbb{R}$ must take the form $s=s_p$.
\end{lem}
\begin{proof}
	Since we assume $V_0=0$, $s$ can be reduced to a morphism $C^{\infty}(U) \rightarrow \mathbb{R}$. Let $x^{\mu}$ be a coordinate system of $U$. Let $f^{\mu}=x^{\mu}-s(x^{\mu})$ and $h=\sum_{\mu} (f^{\mu})^2$. Then $s(h)=0$, which implies that $h$ is non-invertible.	In other words, there exists $p \in U$ such that $x^{\mu}(p)=s(x^{\mu})$ for all $\mu$. Now suppose there exists an $f \in C^{\infty}(U)$ such that $s(f)\neq s_p(f)=f(p)$. Consider the function $h'=h+(f-s(f))^2$. Since $h>0$ for all points of $U/\{p\}$. We know $h'>0$ on $U$. But this contradicts the fact that $s(h')=0$. Hence $s$ must equal $s_p$.
\end{proof}
\begin{thm}\label{maingradom}
	Let $\varphi=(\tilde{\varphi},\varphi^*)$ be a morphism from $\mathcal{U}_1=(U_1,\mathcal{O}_1)$ to $\mathcal{U}_2=(U_2,\mathcal{O}_2)$. Let $(x^{\mu},\theta_{i,a})$ be a coordinate system of $\mathcal{U}_2$. Then $\varphi^*$ is uniquely determined by the equations
	\begin{align*}
		\varphi^*x^{\mu} = y^{\mu},\quad \varphi^*\theta_{i,a} = \eta_{i,a},
	\end{align*}
	where $y^{\mu} \in \mathcal{O}(U_1)_0$, $\eta_{i,a}\in \mathcal{O}(U_1)_i$ and $(\epsilon(y^{\mu}))(p) \in U_2$ for all $p \in U_1$.
\end{thm}
\begin{proof}
	Let $f \in \mathcal{O}_2(U_2)$. By (\ref{decomp}), to construct $\varphi^*f$, we only need to define $\varphi^*f_{\mathcal{J}, \beta}$. But this is straightforward: one just replaces $x^{\mu}$ with $y^\mu$ and $\theta_{i,a}$ with $\eta_{i,a}$ in (\ref{grasscont}). By construction, we have $\varphi^*1 = 1$, $\varphi^*(f+g)=\varphi^*f + \varphi^*g$, and $\varphi^*(fg)=\varphi^*f\varphi^*g$, hence $\varphi^*$ is well-defined. 
	
	Now suppose there exists another $\varphi'^*$ which equals $\varphi^*$ on coordinates. Then they also equals on all polynomials of $(x^{\mu},\theta_{i,a})$. By Lemma \ref{hadamard} and Lemma \ref{critequal}, $\varphi'^* = \varphi^*$.
\end{proof}
\begin{rmk}
	Theorem \ref{maingradom} can be seen as a generalization of the Global Chart Theorem in the $\mathbb{Z}_2$-graded setting (see Theorem 4.2.5 in \cite{Carmeli11}). 
\end{rmk}
\begin{cor}\label{equicat}
	Let $\varphi^*: \mathcal{O}_2(U_2) \rightarrow \mathcal{O}_1(U_1)$ be a ring homomorphism which preserves the $\mathcal{I}$-grading. Then there exists a unique morphism $\varphi': \mathcal{U}_1 \rightarrow \mathcal{U}_2$ such that $\varphi'^* = \varphi^*$.
\end{cor}
\begin{proof}
	First, one can easily show that $\varphi^*$ is actually an $\mathbb{R}$-algebra homomorphism using arguments similar to those in Lemma \ref{unicinf}. Choose a point $p \in U_1$, by Lemma \ref{recovtop}, the morphism $s_p \circ \varphi^*$ must take the form $s_{p'}$ for some $p' \in U_2$. It follows that $\varphi^*(I_{p'}) \subset I_p$. Let $(x^{\mu},\theta_{i,a})$ be a coordinate system of $\mathcal{U}_2$, we then have $\varphi^*x^{\mu}-\epsilon(x^{\mu})(p') \in I_p$. Hence $(\epsilon(\varphi^*x^{\mu}))(p) \in U_2$ for all $p \in U_1$. Next, observe that a coordinate system of $U_2$ restricted to any open subset of it gives a coordinate system of that open subset. 
	Now apply Theorem \ref{maingradom} and Corollary \ref{undermap}.
\end{proof}

\section{Monoidally Graded Manifolds}

\begin{defn}\label{Igm}
	Let $M$ be a $n$-dimensional manifold. An $\mathcal{I}$-graded manifold $\mathcal{M}$ of dimension $n|(m_{i})_{i \in \mathcal{I}}$ is an $\mathcal{I}$-graded ringed space $(M,\mathcal{O}_M)$ which is locally isomorphic to an $\mathcal{I}$-graded domain of dimension $n|(m_{i})_{i \in \mathcal{I}}$. That is, for each $x \in M$, there exist an open neighborhood $U_x$ of $x$,  an $\mathcal{I}$-graded domain $\mathcal{U}$, and an isomorphism of locally ringed spaces
	\begin{align*}
		\varphi=(\tilde\varphi, \varphi^*): (U_x, \mathcal{O}_M|_{U_x}) \rightarrow \mathcal{U}.
	\end{align*}
	$\varphi$ is called a chart of $\mathcal{M}$ on $U_x$.\footnote{We often refer to $U_x$ as a chart too.}
\end{defn}
$M$ with the sheaf $C^{\infty}$ of smooth functions on $M$ is an $\mathcal{I}$-graded manifold of dimension $n|(0,\cdots)$, which is denoted again by $M$ for simplicity. We call $M$ together with a morphism $\mathcal{O} \rightarrow C^{\infty}$ an underlying manifold of $\mathcal{M}$. Equivalently, an underlying manifold of $\mathcal{M}$ is a morphism $\varphi: M \rightarrow \mathcal{M}$ with $\tilde{\varphi}=\mathrm{id}$. 

Let $x \in M$. An open neighborhood $U$ of $x$ on which $\mathcal{O}(U) \cong C^{\infty}(U) \otimes \overline{\mathrm{S}(V)}$ is called a splitting neighborhood. Clearly, every chart is a splitting neighborhood, but not vice versa. The set of splitting neighborhoods form a base of the topology of $M$. For a splitting $U$, there exists sub-algebras $C(U)$ and $D(U)$ of $\mathcal{O}(U)$ such that $C(U)\cong C^{\infty}(U)$, $D(U) \cong \overline{\mathrm{S}(V)}$ and $\mathcal{O}(U)=C(U)\otimes D(U)$. This induces an epimorphism
\begin{align*}
	\epsilon: \mathcal{O}(U) \rightarrow C^{\infty}(U)
\end{align*}
of graded commutative $\mathbb{R}$-algebras, which is a body map of $\mathcal{O}(U)$.

\begin{defn}
	A local coordinate system of $\mathcal{M}$ is the data $(U,x^{\mu},\theta_{i,a})$ where
	\begin{enumerate}
		\item $U$ is a splitting neighborhood of $\mathcal{M}$;
		\item $x^1,\dots,x^n$ are elements of $C(U)$ such that $\epsilon(x^1),\dots,\epsilon(x^n)$ are local coordinate functions of $M$ on $U$;
		\item $\theta_{i,a}$ are homogeneous elements of $\mathcal{D}(U)$ of degree $d(\theta_{i,a})=i$, $i \neq 0$ and $a=1,\cdots,m_i$, which generate $\mathcal{O}(U)$ as a $C(U)$-algebra.
	\end{enumerate}
\end{defn}

\begin{rmk}
	By Theorem \ref{maingradom}, every local coordinate system determines (non-canonically) a chart.
\end{rmk}

Now, let $U$ be an arbitrary open subset of $M$. We can choose a collection of charts $\{U_{\alpha}\}$ such that $U=\bigcup_{\alpha} U_{\alpha}$. For $f \in \mathcal{O}(U)$, one can apply the restriction morphisms to $f$ to get a sequence of sections $f_{\alpha}$ in $\mathcal{O}(U_{\alpha})$. Now, apply $\epsilon$ to each of them to get a sequence of smooth functions $\tilde{f}_{\alpha}$ in $C^{\infty}(U_{\alpha})$. By Proposition \ref{uniepshf}, $\tilde{f}_{\alpha}$ must be compatible with each other, hence can be glued together to get a smooth function $\tilde{f}$ over $U$. In this way, we construct a body map for every open subset of $M$, which are compatible with restrictions. In other words, $\epsilon$ can be seen as a sheaf morphism from $\mathcal{O}$ to $C^{\infty}$.
\begin{prop}\label{uniepsman}
	There exists a unique monomorphism $\varphi: M \rightarrow \mathcal{M}$ with $\tilde{\varphi}=\mathrm{id}$. 
\end{prop}
\begin{proof}
	Existence is guaranteed by $\epsilon$. Uniqueness follows from Proposition \ref{uniepshf}.
\end{proof}
\begin{prop}\label{compepsi}
	Let $\varphi=(\tilde{\varphi}, \varphi^*)$ be a morphism from $\mathcal{M}=(M, \mathcal{O}_M)$ to $\mathcal{N}=(N, \mathcal{O}_N)$. The following diagram commutes.
	\[
	\begin{tikzcd}
		\mathcal{M} \arrow{r}{\varphi}  &\mathcal{N}  \\
		M \arrow[hook]{u} \arrow{r}{\tilde{\varphi}} &N \arrow[hook]{u}
	\end{tikzcd}
	\]
\end{prop}
\begin{proof}
	The proof is essentially the same as the one of Proposition \ref{epicom}.
\end{proof}
\begin{lem}\label{compgramfd}
	Let $\mathcal{O}^1$ be the kernel of $\epsilon$. $\mathcal{O}$ is $\mathcal{O}^1$-adic complete. That is, for any open subset $U$, $\mathcal{O}(U)$ is $\mathcal{O}^1(U)$-adic complete.
\end{lem}
\begin{proof}
	Let $\widehat{\mathcal{O}}$ be the $\mathcal{O}^1$-adic completion of $\mathcal{O}$.\footnote{For each open subset $U$, one has $\widehat{\mathcal{O}}(U) = \varprojlim \mathcal{O}(U)/\mathcal{O}^n(U)$, where $\mathcal{O}^n(U)$ is the $n$-th power of $\mathcal{O}^1(U)$.} There exists a canonical morphism $\iota: \mathcal{O} \rightarrow \widehat{\mathcal{O}}$. Since $\mathcal{O}$ is locally $\mathcal{O}^1$-adic complete, the induced stalk morphism $\iota_p: \mathcal{O}_p \rightarrow \widehat{\mathcal{O}}_p$ is an isomorphism for each $p \in M$. It follows that $\mathcal{O}$ is $\mathcal{O}^1$-adic complete.
\end{proof}
\begin{defn}
	An $\mathcal{I}$-graded manifold $\mathcal{M}$ is called projected if there exists a splitting of the short exact sequence of sheaves of rings
	\begin{align}\label{proj}
		0 \longrightarrow \mathcal{O}^1 \longrightarrow \mathcal{O} \xrightarrow{~\epsilon~} C^{\infty} \longrightarrow 0,
	\end{align}
	where $\mathcal{O}^1$ is the kernel of $\epsilon$.
\end{defn}
The structure sheaf $\mathcal{O}$ of a projected manifold is a $C^{\infty}$-module.
\begin{defn}
	A projected $\mathcal{I}$-graded manifold $\mathcal{M}$ is called split if there exists a splitting of the short exact sequence of $C^{\infty}$-modules
	\begin{align}\label{split}
		0 \longrightarrow \mathcal{O}^2 \longrightarrow \mathcal{O}^1 \xrightarrow{~\pi~} \mathcal{O}^1/\mathcal{O}^2  \longrightarrow 0,
	\end{align}
	where $\mathcal{O}^2$ is the square of $\mathcal{O}^1$, $\pi$ is the canonical quotient map.
\end{defn}
Let $\mathcal{O}$ be the structure sheaf of a projected $\mathcal{I}$-graded manifold. Let $\mathcal{F}$ denote the sheaf $\mathcal{O}^1/\mathcal{O}^2$. $\mathcal{F}$ is an $\mathcal{I}$-graded $C^{\infty}$-module and we can define its formal symmetric power $\overline{\mathrm{S}(\mathcal{F})}$. By construction, the ringed space  $\mathcal{M}_{S}=(M,\overline{\mathrm{S}(\mathcal{F})})$ is also a projected $\mathcal{I}$-graded manifold.
\begin{lem}
	Let $\mathcal{M}=(M,\mathcal{O})$ be a projected $\mathcal{I}$-graded manifold. $\mathcal{M}$ is split if and only if $\mathcal{M} \cong \mathcal{M}_{S}$.
\end{lem}
\begin{proof}
	Let $\iota: \mathcal{F} \rightarrow \overline{\mathrm{S}(\mathcal{F})}$ be the canonical monomorphism. $\iota$ splits the short exact sequence (\ref{split}). Now if $\mathcal{M}$ is split, one can find a monomorphism $F: \mathcal{F} \rightarrow \mathcal{O}$ of $C^{\infty}$-modules such that $F(\mathcal{F}(U)) \subset \mathcal{O}^1(U)$ for any open subset $U$. By Lemma \ref{uniforshf} and Lemma \ref{compgramfd}, there exists a unique $C^{\infty}$-algebra morphism $\tilde{F}: \overline{\mathrm{S}(\mathcal{F})} \rightarrow \mathcal{O}$ such that $\tilde{F} \circ \iota = F$. By Remark \ref{isofor}, $\tilde{F}$ induces an isomorphism for each stalk. Hence $\mathcal{M} \cong \mathcal{M}_{S}$.
\end{proof}
\begin{lem}
	Every projected $\mathcal{I}$-graded manifold is split.
\end{lem}
\begin{proof}
	Due to the existence of a smooth partition of unity on $M$, $H^q(M, \mathrm{Hom}(\mathcal{O}^1/\mathcal{O}^2, \mathcal{O}^2))$ vanishes for $q \geq 1$. By Lemma \ref{obspl}, there is no obstruction of the existence of a splitting of (\ref{split}).
\end{proof}
\begin{lem}
	Every $\mathcal{I}$-graded manifold is projected.
\end{lem}
\begin{proof}
	Let $\mathcal{O}_{(i)} = \mathcal{O}/\mathcal{O}^{i+1}$. Let $\phi_{(0)}: C^{\infty} \rightarrow \mathcal{O}_{(0)}$ be the identity. (By Proposition \ref{uniepsman}, there is a unique identification $\mathcal{O}_{(0)} \cong C^{\infty}$.) One can construct by induction on $i$ mappings $\phi_{(i+1)}: C^{\infty} \rightarrow \mathcal{O}_{(i+1)}$ such that $\pi_{i+1,i} \circ \phi_{(i+1)} = \phi_{(i)}$, where $\pi_{i+1,i}: \mathcal{O}_{i+1} \rightarrow \mathcal{O}_i$ is the canonical epimorphism. As is shown in \cite{Manin97}, one can construct an element
	\begin{align*}
		\omega(\phi_{(i)}) \in H^1(M, (\mathcal{T} \otimes \mathrm{S}^{i+1}(\mathcal{F}))_0)
	\end{align*}
	as the obstruction to the existence of $\phi_{(i+1)}$, where $\mathcal{T}$ is the tangent sheaf of $M$. Due to the existence of a smooth partition of unity on $M$, $H^1(M, (\mathcal{T} \otimes \mathrm{S}^{i+1}(\mathcal{F}))_0)=0$ and $\omega(\phi_{(i)})=0$. It follows that there exists a unique morphism
	$\phi: C^{\infty} \rightarrow \varprojlim \mathcal{O}_{(i)}$ such that $\pi_i \circ \phi = \phi_{(i)}$, where $\pi_i: \varprojlim \mathcal{O}_{(i)} \rightarrow \mathcal{O}_i$ is the canonical epimorphism. By Lemma \ref{compgramfd}, $\phi$ can be seen as a morphism from $C^{\infty}$ to $\mathcal{O}$. Note that $\pi_0 = \epsilon$ and $\pi_0 \circ \phi = \phi_{(0)} = \mathrm{id}$. $\phi$ splits (\ref{proj}).   
\end{proof}
\begin{cor}\label{splcor}
	Every $\mathcal{I}$-graded manifold is split.
\end{cor}
Let $V$ be a (finite dimensional) $\mathcal{I}$-graded vector space. An $\mathcal{I}$-graded vector bundle $\pi: E \rightarrow M$ is a vector bundle such that the local trivialization map
$
\varphi_U: \pi^{-1}(U) \rightarrow U \times V
$
is a morphism of $\mathcal{I}$-graded vector spaces when restricted to $\pi^{-1}(x)$, $x \in U \subset M$. In other words, $E = \bigoplus_{k \in \mathcal{I}} E_k$ where $E_k$ are vector bundles whose fibers consist of elements of degree $k$. To any $\mathcal{I}$-graded vector bundle $E$ we can associate an $\mathcal{I}$-graded ringed space with the underlying topological space being $M$ and the structure sheaf being the sheaf of sections of $\overline{\mathrm{S}(\bigoplus_{k \in \mathcal{I}} (E_{k})^*)}$. (This is an $\mathcal{I}$-graded manifold in our sense if the fiber of $E$ does not contain elements of degree $0$.) Corollary \ref{splcor} can then be rephrased as
\begin{thm}
	Every $\mathcal{I}$-graded manifold can be obtained from an $\mathcal{I}$-graded vector bundle.
\end{thm}

\section{Vector Fields and Tangent Sheaves}

Throughout this section, every algebra is assumed to be real. 

Let $R$ be a unital associative commutative $\mathcal{I}$-graded algebra. An $\mathcal{I}$-graded algebra $A$ over $R$ is defined to be an $\mathcal{I}$-graded algebra $A$ equipped with a left $R$-module structure such that $R_iA_j \subset A_{i+j}$, and
\begin{align*}
	r(ab)=(ra)b=(-1)^{p(r)p(a)}a(rb)
\end{align*}
for $r\in R$ and $a,b \in A$. Recall that when we write $p(r)p(a)$, we mean $p(d(r)d(a))$. We also require that $1 a = a$, where $1 \in R$ is the identity element.
\begin{defn}
	An $\mathcal{I}$-graded Lie algebra over $R$ is an $\mathcal{I}$-graded algebra $L$ over $R$ whose multiplications (denoted by $[\cdot,\cdot]$) satisfy
	\begin{align}
		&[a,b]=-(-1)^{p(a)p(b)}[b,a], \label{bl1}\\
		&[a,[b,c]]=[[a,b],c] + (-1)^{p(a)p(b)}[b,[a,c]], \label{bl2}
	\end{align} 
	for all $a,b,c \in L$. 
\end{defn}

The space of endomorphisms $\mathrm{Hom}(A,A)$ (or $\mathrm{gl}(A)$) of an $\mathcal{I}$-graded $R$-module $A$ is an associative $K(\mathcal{I})$-graded algebra over $R$. It can be also viewed as a $K(\mathcal{I})$-graded Lie algebra over $R$ by setting
\[
[f,g]=f\circ g -(-1)^{p(f)p(g)}g \circ f
\]
for all $f,g\in \mathrm{Hom}(A,A)$. In the case of $A$ being an $\mathcal{I}$-graded algebra, an endomorphism $D$ is said to be a derivation if 
\begin{align}\label{supder}
	D(ab)=D(a)b+(-1)^{p(D)p(a)}aD(b).
\end{align}
It is easy to check that derivations of $A$ form a $K(\mathcal{I})$-graded Lie subalgebra of $\mathrm{gl}(A)$ over $R$.

\begin{defn}
	Let $\mathcal{M}=(M,\mathcal{O})$ be an $\mathcal{I}$-graded manifold. A (local) vector field over $\mathcal{M}$ is a derivation of $\mathcal{O}(U)$, where $U$ is an open subset of $M$.
\end{defn}

Local vector fields over $\mathcal{M}$ actually form a ($K(\mathcal{I})$-graded) sheaf on $M$. To prove this, we need the partition of unity lemma in the $\mathcal{I}$-graded setting.

\begin{lem}
	Let $f \in \mathcal{O}(M)$ such that $\epsilon(f)(x) \neq 0$ for all $x \in M$. $f$ is invertible.
\end{lem}
\begin{proof}
	Choose an open cover $\{U_{\alpha}\}$ of charts of $M$. Let $f_{\alpha}$ denote $\rho_{U_{\alpha},M}(f)$. Each $f_{\alpha}$ is invertible by Lemma \ref{invfor}. Let $f^{-1}_{\alpha}$ denote the inverse of  $f_{\alpha}$. By uniqueness of the inverse, $f^{-1}_{\alpha}$ are compatible with each other, hence can be glued to give a section $f^{-1} \in \mathcal{O}(M)$, which is the inverse of $f$.
\end{proof}
\begin{lem}\label{pou}
	Let $\{U_{\alpha}\}$ be an open cover of $M$. There exists a locally finite refinement $\{V_{\beta}\}$ of $\{U_{\alpha}\}$ and a family of functions $\{l_{\beta} \in \mathcal{O}(M)_0\}$ such that 
	\begin{enumerate}
		\item $\mathrm{supp}~l_{\beta} \subset V_{\beta}$ is compact and $\epsilon(l_{\beta}) \geq 0$ for all $\beta$;
		\item $\sum_{\beta} l_{\beta} = 1$.
	\end{enumerate} 
\end{lem}
\begin{proof}
	First, find a partition of unity $\{\tilde{l}_{\beta}\}$ of $M$ subordinate to $\{V_{\beta}\}$. Choose $l'_{\beta} \in \mathcal{O}(V_{\beta})$ such that $\epsilon(l'_{\beta})=\tilde{l}_{\beta}$. Since $\tilde{l}_{\beta}$ are invertible, we can then set $l_\beta$ to be $(\sum_{\beta} l'_{\beta})^{-1}l'_{\beta}$.
\end{proof}
Using Lemma \ref{pou}, it is not hard to prove the following lemma.
\begin{lem}\label{localivec}
	Let $U$ and $V$ be open in $M$ such that $V \subset U$. Let $D$ be a derivation of $\mathcal{O}(U)$. Then there exists a unique derivation $D'$ of $\mathcal{O}(V)$ such that $D'(\rho_{V,U}(f))=\rho_{V,U}(D(f))$ for all $f \in \mathcal{O}(U)$.
\end{lem}
We skip the proof of Lemma \ref{localivec} since it is essentially the same as the one in the $\mathbb{Z}_2$-graded setting \cite{Leites80}.  Note that Lemma \ref{localivec} implies that local vector fields over $\mathcal{M}$ form a presheaf $\mathfrak{X}$ on $M$.
\begin{prop}
	$\mathfrak{X}$ is a sheaf on $M$. 
\end{prop} 
\begin{rmk}
	$\mathfrak{X}$ is called the tangent sheaf of $\mathcal{M}$.
\end{rmk}
\begin{proof}
	Let $U$ be an open subset of $M$ with an open cover $\{U_{\alpha}\}$. Let $D_{\alpha} \in \mathfrak{X}(U_{\alpha})$ be compatible with each other. We obtain a $D \in \mathfrak{X}(U)$ by setting $D(f)$ to be unique function obtained by gluing $D_{\alpha}(f_{\alpha})$, where $f \in \mathcal{O}(U)$ and $f_{
		\alpha}=\rho_{U_{\alpha},U}(f)$. $D(f)$ is well defined because $\rho_{U_{\alpha}\cap U_{\beta}, U_{\alpha}}(D_{\alpha}(f_{\alpha})) = \rho_{U_{\alpha}\cap U_{\beta}, U_{\alpha}}(D_{\alpha})(\rho_{U_{\alpha}\cap U_{\beta}, U_{\alpha}}(f_{\alpha}))=\rho_{U_{\alpha}\cap U_{\beta}, U_{\beta}}(D_{\beta})(\rho_{U_{\alpha}\cap U_{\beta}, U_{\beta}}(f_{\beta}))=\rho_{U_{\alpha}\cap U_{\beta}, U_{\beta}}(D_{\beta}(f_{\beta}))$.
\end{proof}
We are now ready to give the following definition.
\begin{defn}
	A $QK$-manifold is a bigraded manifold equipped with three vector fields $Q$, $K$ and $d$ of degree $(0,1)$, $(1,-1)$ and $(1,0)$, respectively, satisfying the following relations
	\[
	Q^2=0, \quad QK+KQ=d, \quad Kd+dK=0.
	\]
\end{defn}
$QK$-manifolds can be used to study the descent equations \eqref{desceq} in a cohomological field theory. The physical observables $\mathcal{O}^{(p)}$ can be interpreted as functions over a $QK$-manifold. Fix a $Q$-closed function $\mathcal{O}^{(0)}$, a $K$-sequence is defined by setting $\mathcal{O}^{(p)}=\frac{1}{p!}K^p \mathcal{O}^{(0)}$. A sequence $\{\mathcal{O}^{(p)}\}_{p=0}^n$ is called an exact sequence if there exists another sequence $\{\mathcal{P}^{(p)}\}_{p=0}^n$ such that $\mathcal{O}^{(p)}=Q \mathcal{P}^{(p)} + d \mathcal{P}^{(p-1)}$ for $p \geq 1$ and $\mathcal{O}^{(0)}=Q \mathcal{P}^{(0)}$. In \cite{Jiang22}, the author proved that
\begin{thm}
	Every solution to \eqref{desceq} is a $K$-sequence up to an exact sequence.
\end{thm}

\section{Conclusions}

In this paper, we have given a definition of $\mathcal{I}$-graded manifolds, where $\mathcal{I}$ is a countable cancellative commutative semi-ring. Such a definition unifies different objects such as supermanifolds, graded manifolds and colored supermanifolds. We have proved the existence and uniqueness of an underlying manifold of an $\mathcal{I}$-graded manifold. Furthermore, we have also proved Batchelor's theorem in this generalized setting, namely that every $\mathcal{I}$-graded manifold can be obtained from an $\mathcal{I}$-graded vector bundle. At the end of this paper, we have discussed a special class of bigraded manifolds, the $QK$-manifolds, and their applications to cohomological field theories.


\section*{Acknowledgement}

The author would like to thank Enno Keßler and Sylvain Lavau for pointing out to him that Batchelor's theorem was proved in both the $\mathbb{Z}_2^n$-graded setting and the $\mathbb{Z}$-graded setting. He would also like to thank an anonymous referee for pointing out the connection between \eqref{grasscont} and the definition of a supersmooth function in the sense of Rogers-DeWitt. This work was supported by the International Max Planck Research School Mathematics in the Sciences.


\begin{bibsection}
	\begin{biblist}
		
	\bib{Kostant77}{incollection}{
		title={Graded manifolds, graded Lie theory, and prequantization},
		author={Kostant, B.},
		booktitle={Differential Geometrical Methods in Mathematical Physics},
		pages={177--306},
		date={1977},
		publisher={Springer Berlin, Heidelberg}
	}
    \bib{Witten82}{article}{
    	title={Supersymmetry and Morse theory},
    	author={Witten, E.},
    	journal={Journal of Differential Geometry},
    	volume={17},
    	number={4},
    	pages={661--692},
    	date={1982},
    	publisher={Lehigh University}
    }
    \bib{Cattaneo11}{article}{
    	author={Cattaneo, A. S.},
    	author={Sch\"{a}tz, F.},
    	title={Introduction to supergeometry},
    	journal={Reviews in Mathematical Physics},
    	volume={23},
    	number={06},
    	pages={669-690},
    	date={2011}    
    }
	\bib{Fairon17}{article}{
		title={Introduction to graded geometry},
		author={Fairon, M.},
		journal={European Journal of Mathematics},
		volume={3},
		number={2},
		pages={208--222},
		date={2017},
		publisher={Springer}
	}	
    \bib{Alexandrov1997}{article}{
    	title={The geometry of the master equation and topological quantum field theory},
    	author={Alexandrov, M.},
    	author={Schwarz, A.},
    	author={Zaboronsky, O.},
    	author={Kontsevich, M.},
    	journal={International Journal of Modern Physics A},
    	volume={12},
    	number={07},
    	pages={1405--1429},
    	date={1997},
    	publisher={World Scientific}
    }
    \bib{Witten88b}{article}{
    	title={Topological sigma models},
    	author={Witten, E.},
    	journal={Communications in Mathematical Physics},
    	volume={118},
    	number={3},
    	pages={411--449},
    	date={1988},
    	publisher={Springer}
    }
    \bib{Witten88a}{article}{
    	title={Topological quantum field theory},
    	author={Witten, E.},
    	journal={Communications in Mathematical Physics},
    	volume={117},
    	number={3},
    	pages={353--386},
    	date={1988},
    	publisher={Springer}
    }
    \bib{Jiang22}{article}{
    	title={Mathematical structures of cohomological field theories},
    	author={Jiang, S.},
    	eprint={arxiv:2202.12425},
    	date={2022}
    }
    \bib{Leites80}{article}{
    	title={Introduction to the theory of supermanifolds},
    	author={Leites, D. A.},
    	journal={Russian Mathematical Surveys},
    	volume={35},
    	number={1},
    	pages={1--64},
    	date={1980},
    	publisher={IOP Publishing}
    }
    \bib{Manin97}{book}{
    	title={Gauge Field Theory and Complex Geometry},
    	author={Manin, Y. I.},
    	volume={289},
    	pages={},
    	date={1997},
    	publisher={Springer Berlin, Heidelberg}
    }
    \bib{Carmeli11}{book}{
    	title={Mathematical Foundations of Supersymmetry},
    	author={Carmeli, C.},
    	author={Caston, L.},
    	author={Fioresi, R.},
    	volume={15},
    	date={2011},
    	publisher={European Mathematical Society}
    }
    \bib{Bartocci12}{book}{
    	title={The Geometry of Supermanifolds},
    	author={Bartocci, C.},
    	author={Bruzzo, U.},
    	author={Ruip{\'e}rez, D. H.},
    	date={2012},
    	publisher={Springer Dordrecht}
    }
    \bib{Kessler19}{book}{
    	author = {Ke{\ss}ler, E.},
    	title = {Supergeometry, Super Riemann Surfaces and the Superconformal Action Functional},
    	publisher = {Springer Cham},
    	volume = {2230},
    	date = {2019}
    }
    \bib{Batchelor79}{article}{
    	title={The structure of supermanifolds},
    	author={Batchelor, M.},
    	journal={Transactions of the American Mathematical Society},
    	volume={253},
    	pages={329--338},
    	date={1979}
    }
    \bib{Covolo16}{article}{
    	title={Splitting theorem for $\mathbb{Z}_2^n$-supermanifolds},
    	author={Covolo, T.},
    	author={Grabowski, J.},
    	author={Poncin, N.},
    	journal={Journal of Geometry and Physics},
    	volume={110},
    	pages={393--401},
    	date={2016},
    	publisher={Elsevier}
    }
    \bib{Kotov21}{article}{
    	title={The category of $\mathbb{Z}$-graded manifolds: what happens if you do not stay positive},
    	author={Kotov, A.},
    	author={Salnikov, V.},
    	eprint={arxiv:2108.13496},
    	date={2021}
    }
    \bib{Singh11}{book}{
    	title={Basic Commutative Algebra},
    	author={Singh, B.},
    	date={2011},
    	publisher={World Scientific}
    }
    \bib{Rogers07}{book}{
    	title={Supermanifolds: Theory and Applications},
    	author={Rogers, A.},
    	date={2007},
    	publisher={World Scientific}
    }
	\end{biblist}
\end{bibsection}
\end{document}